\documentclass{article}

\usepackage{delarray,verbatim,enumerate,a4wide}
\usepackage{amsmath,amsthm,amstext,amsbsy,amssymb,amsfonts,amscd}
\usepackage[all]{xy}
\usepackage[french]{babel}
\usepackage[T1]{fontenc}
\usepackage[utf8]{inputenc}

\usepackage{scalerel}

\usepackage{upgreek,bbm}

\usepackage[small]{titlesec}

\usepackage{color}
\definecolor{vert}{rgb}{0.1,0.4,0.2}
\usepackage[colorlinks=true,linkcolor=blue,citecolor=vert]{hyperref}

\usepackage{calligra}
\DeclareFontShape{T1}{calligra}{m}{n}{<->s*[0.95]callig15}{}
\DeclareMathAlphabet{\mathscr}{T1}{calligra}{m}{n}

\newtheorem{Th}{Théorème}

\newtheorem{Def} [Th]{Définition}

\newtheorem*{TP}{Théorème principal}
\newtheorem*{Sco*}{Scolie}

\def\Remarque{\smallskip\noindent {\it Remarque.~}}

\def\Nota{\smallskip\noindent {\it Nota.~}}

	\def\ZZ{\mathbb Z}		
\def\F2{\mathbb{F}_2}	\def\Z2{\mathbb{Z}_2}		
\def\Zl{{\mathbb{Z}_\ell}} 			

 				\def\U{\mathcal  U}	\def\F{\mathcal  F}	
\def\J{\mathcal  J}  				\def\R{\mathcal  R}		\def\Z{\mathcal  Z}
\def\Dl{\mathcal  D\ell} 	  	\def\Cl{\mathcal  C\!\ell}	
				\def\Dl{\mathcal  D\ell}

		\def\p{{\mathfrak p}}						
		\def\l{{\mathfrak l}}		\def\d{{\mathfrak d}}

\def\wi{\widetilde}				
	
	\def\deg{\operatorname{deg}}		
\def\Gal{\operatorname{Gal}}			
\def\Ker{\operatorname{Ker}}				

\newcommand\scale[2]{\vstretch{#1}{\hstretch{#1}{#2}}}

\newcommand\si[1]{\scale{.7}{#1}}	
\newcommand\ph{{\phantom{*}}}
	\newcommand\lc{{\scale{.8}{\rm lc}}}	\newcommand\nr{{\scale{.8}{\rm nr}}}

\newcommand\res{{\scale{.8}{\rm res}}}	\newcommand\tor{{\scale{.8}{\rm tor}}}
\def\%{{\scale{.8}{\infty}}}

\makeatletter
\newcommand*\wt[2][0.2ex]{%
        \begingroup
        \mathchoice{\wt@helper{#1}{#2}{\displaystyle}{\textfont}}
                   {\wt@helper{#1}{#2}{\textstyle}{\textfont}}
                   {\wt@helper{#1}{#2}{\scriptstyle}{\scriptfont}}
                   {\wt@helper{#1}{#2}{\scriptscriptstyle}{\scriptscriptfont}}%
        \endgroup
        #2%
}
\newcommand*\wt@helper[4]{%
        \def\currentfont{\the#41}%
        \def\currentskewchar{\char\the\skewchar\currentfont}%
        \setbox\tw@\hbox{\currentfont$#2$\currentskewchar}%
        \dimen@ii\wd\tw@
        \setbox\tw@\hbox{\currentfont$#2${}\currentskewchar}%
        \advance\dimen@ii-\wd\tw@
        \rlap{\raisebox{-#1}{$\m@th#3\kern\dimen@ii\widetilde{\phantom{#2}}$}}%
}
\makeatother

\def\wE{\,\wt[0.1ex]{\!\mathcal E}}		\def\wU{\wt[0.2ex]{\mathcal U}}
\def\wJ{\,\wt[0.2ex]{\!\mathcal J}}	\def\wCl{\wt[0.1ex]{\mathcal C\!\ell}} 
				
\def\wnu{\wt[0.1ex]{\nu}}			\def\we{\wt[0.1ex]{e}}
			
			\def\wH{\,\wt[0.2ex]{H}}			
		
\date{}

\title{Classes logarithmiques et capitulation II}

\author{Jean-François {\sc Jaulent}}

\begin{document}
\maketitle
\medskip

{\footnotesize
\noindent{\bf Résumé.} Nous établissons un analogue du Théorème d'Artin-Furtwängler sur la capitulation en transposant aux classes logarithmiques  la preuve algébrique classique du Théorème de l'idéal principal.\smallskip

\noindent{\bf Abstract.} We establish a logarithmic version of the classical result of Artin-Furwängler on the principalization of ideal classes in the Hilbert class-field by applying the group theoretic description of the transfert map.\smallskip


}


\section{Introduction}

L'objet de cette note est d'établir pour les groupes de classes logarithmiques de degré nul $\,\wCl_K$ des corps de nombres un analogue du Théorème d'Artin-Furwängler de 1930 (cf. e.g. \cite{Art,A-T,Fur}), lequel affirme que le groupe des classes d'idéaux $Cl_K$ d'un tel corps $K$ capitule dans son corps des classes de Hilbert $L=H_K$; en d'autres termes que les idéaux de $K$ se principalisent dans son extension abélienne non ramifiée maximale $L$.\smallskip

Ce résultat célèbre, qui relève essentiellement de la Théorie du corps de classes, trouve son origine dans le Théorème 94 de Hilbert \cite{Hil} qui dit que dans une extension cyclique non-ramifiée de degré premier $L/K$ il existe une classe non-triviale de $Cl_K$ qui se capitule dans $Cl_L$.\smallskip

Or, comme exposé dans une étude précédente sur cette question \cite{J54}, ce dernier théorème est en défaut pour les groupes de classes logarithmiques, puisque les calculs de K. Belabas \cite{BJ} effectués à l'aide du logiciel {\sc pari} fournissent des exemples de corps quadratiques $K$  ayant un 3-groupe de classes logarithmiques $\,\wCl_K$ d'ordre 3 mais admettant des extensions cycliques logarithmiquement non-ramifiées $L$ de degré 3 dans lesquelles $\,\wCl_K$ ne capitule pas. Qui plus est, dans cette configuration, $K$ possède plusieurs 3-extensions cycliques qui sont logarithmiquement non-ramifiées, ce qui rend problématique la définition même d'un 3-corps de classes de Hilbert au sens logarithmique.\smallskip

C'est ce double problème de définition et de capitulation que nous nous résolvons ici:

\begin{TP}
Soient $K$ un corps de nombres, $\ell$ un nombre premier et $\ell^{\,\we_{\si{K}}}$ l'exposant du sous-groupe de torsion $\,\wCl^{\,\tor}_K$ du pro-$\ell$-groupe des classes logarithmiques (de degré arbitraire) de $K$. Alors le corps $K$ possède une $\ell$-extension abélienne naturelle $\wH_K$ logarithmiquement non-ramifiée et d'exposant $\ell^{\we_{\si{K}}}\!$, maximale sous ces conditions ; laquelle est telle que $\,\wCl^{\,\tor}_K$ capitule dans $\,\wCl^{\,\tor}_L$.\par

Nous disons que $\wH_K$ est le $\ell$-corps de classes logarithmiques normalisé de $K$.
\end{TP}

\begin{Sco*}
Il suit de là que $\,\wCl^{\,\tor}_K$ capitule dans la pro-$\ell$-extension abélienne $K^\lc$ logarithmiquement non-ramifiée maximale de $K$.
\end{Sco*}

\Remarque Nous avons choisi de donner ci-dessus une formulation inconditionnelle du Théorème principal. Mais il est possible d'être plus précis sous les conjectures $\ell$-adiques standard: comme expliqué plus loin, la conjecture de Gross-Kuz'min \cite{Gro,Kuz} pour un corps $K$ revient à postuler que le sous-groupe $\,\wCl_K$ des classes de degré nul du groupe des classes logarithmiques (de degré arbitraire) de $K$ est précisément son plus grand sous-module fini; autrement dit que l'on a: $\wCl^{\,\tor}_K=\,\wCl_K$.\par

Par ailleurs, la principalisation d'un diviseur $\d_{\si{K}}$ (au sens ordinaire comme logarithmique) dans une extension algébrique arbitraire a toujours lieu de fait dans une sous-extension de degré fini, à savoir celle, disons $L=K[\alpha]$, engendrée par le générateur obtenu. Il en résulte que $\d_{\si{K}}^{\si{[L:K]}}$ est engendré par $N_{\si{L/K}}(\alpha)$ et donc que la classe d'un tel $\d_{\si{K}}$ est bien d'ordre fini.

\newpage

\section{Construction du corps de classes logarithmiques normalisé}

Le pro-$\ell$-groupe des classes logarithmiques $\,\Cl_K$ a été introduit dans \cite{J28} par analogie avec le groupe des classes au sens habituel en  envoyant le tensorisé $\R_K=\Zl\otimes_\ZZ K^\times$ dans le $\Zl$-module des diviseurs $\Dl=\oplus_{\p}\Zl\,\p$ construit sur les places finies de $K$ par la famille $(\wnu_\p)_\p$ obtenue en remplaçant les valuations habituelles $\nu_\l$ aux places $\l$ qui divisent $\ell$ par les valuations $\ell$-adiques $\wnu_\l$ données (à normalisation près) par les logarithmes des valeurs absolues $\ell$-adiques $\wnu_\l(\cdot)= \log_\ell(|\cdot|)/\log_\ell(1+\ell)$:\smallskip

\centerline{$1 \rightarrow \wE_K \rightarrow  \R_K\overset{\wi\nu}{\longrightarrow}  \Dl_K \rightarrow \Cl_K \rightarrow 1$.}\smallskip

\noindent Dans la suite exacte obtenue le noyau $\,\wE_K$ à gauche est ainsi le groupe des {\em unités logarithmiques}; le conoyau $\,\Cl_K$ à droite, celui des {\em classes logarithmiques}. Contrairement aux groupes de classes et d'unités au sens habituel, ce sont donc par construction des objets $\ell$-adiques. Leur importance provient de leur interprétation par la Théorie $\ell$-adique du corps de classes (cf. \cite{Gra,J28,J31}).\smallskip

Pour voir cela, introduisons pour chaque place finie $\p$ de $K$  le compactifié $\ell$-adique du groupe $K_\p^\times$ défini par $\R_\p=\varprojlim K_\p^\times/K_\p^{\times\ell^n}$; et notons $\J_K= \prod_\p^\res \R_\p$ le $\ell$-adifié du groupe des idèles.\smallskip

 Du point de vue local, le noyau $\U_\p$ de $\nu_\p$ dans $\R_\p$ (autrement dit le sous-groupe des unités de $\R_\p$) est le groupe de normes associé à la $\Zl$-extension non ramifiée de $K_\p$; tandis que le noyau $\wU_\p$ de $\wnu_\p$ (i.e. le sous-groupe des unités logarithmiques locales) correspond, lui, à sa $\Zl$-extension cyclotomique. En d'autres termes une $\ell$-extension abélienne de $K_\p$ est {\em logarithmiquement non-ramifiée} si et seulement si elle est contenue dans la $\Zl$-extension cyclotomique $K^c_\p$ de $K_\p$.\smallskip
 
Du point de vue global, le $\ell$-groupe des classes d'idéaux s'interprète comme groupe de Galois de la $\ell$-extension abélienne non ramifiée maximale $K^{\nr}$ de $K$; et le $\ell$-groupe des classes logarithmiques $\,\Cl_K\simeq \J_K/\prod_\p\wU_\p\R_K$ comme groupe de Galois de sa pro-$\ell$-extension abélienne localement cyclotomique maximale $K^\lc$. Le corps $K^\lc$ est ainsi la plus grande pro-$\ell$-extension abélienne de $K$ qui est complètement décomposée (en toutes ses places)  au-dessus de la $\Zl$-extension cyclotomique $K_\%$. Et comme $K^\lc$ contient $K_\%$, le groupe  $\,\Cl_K$ n'est donc jamais fini.\smallskip

La surjection de $\J_K$ dans le groupe $\Gamma_{\!K}=\Gal(K_\%/K) \simeq \Zl$ fournit alors un morphisme {\em degré}:\smallskip

\centerline{$\deg_K:\J_K\twoheadrightarrow \Zl$;}

\noindent dont le noyau $\wJ_K$ est, par construction, le sous-groupe normique de $\J_K$ attaché à $K_\%$. Le quotient\smallskip

\centerline{$\,\wCl_K \simeq \wJ_K/\prod_\p\wU_\p\R_K \simeq \Gal(K^\lc/K_\%)$}\smallskip

\noindent est ainsi le {\em sous-groupe des classes logarithmiques de degré nul}, qui est l'objet de cette note.\smallskip

En particulier, le sous-module de torsion $\,\wCl_K^{\,\tor}$ de $\,\Cl_K$ est contenu dans $\,\wCl_K$.
Et la {\em Conjecture de Gross-Kuz'min} (pour le corps $K$ et le premier $\ell$), qui revient à postuler la finitude de $\,\wCl_K$, affirme l'égalité $\,\wCl_K^{\,\tor}=\wCl_K$. Vérifiée en particulier lorsque $K$ est abélien en vertu du résultat de Baker-Brumer (cf. e.g. \cite{Gra,Gre,J28}), c'est, comme expliqué dans \cite{J10}, une conséquence d'une conjecture plus générale d'indépendance $\ell$-adique de nombres algébriques, qui résulte elle-même de la conjecture de Schanuel $\ell$-adique.\smallskip

Inconditionnellement, en revanche, nous avons donc:\smallskip

\centerline{$\Cl_K \simeq \Zl^{1+\delta_{\si{K}}}\! \oplus \,\wCl^{\,\tor}_K$,}\smallskip

\noindent où $\delta_{\si{K}}$ mesure le défaut de la conjecture de Gross-Kuz'min dans $K$ relativement au premier $\ell$.

\begin{Def}
Étant donnés un corps de nombres $K$ et un nombre premier $\ell$, notons $\ell^{\we_{\si{K}}}\!$ l'exposant du sous-groupe de torsion $\,\Cl_K^{\,\tor}$ du pro-$\ell$-groupe  $\,\Cl_K$ des classes logarithmiques (de degré arbitraire).\smallskip

\noindent Nous disons que l'extension abélienne d'exposant $\ell^{\we_{\si{K}}}$ logarithmiquement non-ramifiée $\wH_K$ fixée par le sous-groupe $\J_K^{\ell^{\we_{\si{K}}}}\wU_K \R_K$ de $\J_K$ est le $\ell$-corps de classes logarithmiques normalisé de $K$.\smallskip

Le corps $\wH_K$ est donc, par construction, la plus grande $\ell$-extension abélienne de $K$ qui est d'exposant $\ell^{\we_{\si{K}}}\!$\! et logarithmiquement non-ramifiée.
\end{Def}

\Nota La conjecture de Gross-Kuz'min pour $K$ postule de façon équivalente $\delta_{\si{K}}=0\;$ et $\,\wCl_K^{\,\tor}=\wCl_K$.

\section{Preuve sous la conjecture de Gross-Kuz'min}

Partons d'un corps de nombre arbitraire $K$, notons $L=\wH_K$ son $\ell$-corps de classes logarithmiques normalisé, puis $M=\wH_L$ celui de $L$.\smallskip

Faisons l'hypothèse dans cette section que le corps $L$ satisfait la conjecture de Gross-Kuz'min (pour le premier $\ell$), autrement dit que $L$ possède une unique $\Zl$-extension localement cyclotomique, à savoir sa $\Zl$-extension cyclotomique $L_\%$, de sorte qu'il en est de même pour son sous-corps $K$. 

Écrivons $\Gamma_{\!K}=\Gal(K_\%/K)$; de même $\Gamma_{\!L}=\Gal(L_\%/L)$; et, pour chaque entier naturel $n$, notons $K_n$ l'unique sous-extension de degré $\ell^n$ de la $\Zl$-extension cyclotomique $K_\%$ de $K$ et $L_n$ celle de $L_\%$. Soit enfin $\ell^{\we_{\si{K}}}$ l'exposant de $\,\wCl_K\simeq\Gal(K^\lc/K_\%)$ et $\ell^{\we_{\si{L}}}$ celui de $\,\wCl_L\simeq\Gal(L^\lc/L_\%)$.\smallskip

Par construction, le corps $M=\wH_L$ est une extension galoisienne logarithmiquement non-ramifiée de $K$. Soit donc $G=\Gal(M/K)$ son groupe de Galois et $G'$ le groupe dérivé. Le sous-corps des points fixes de $G'$ est ainsi la sous-extension maximale de $K^\lc=L_\%$ qui est d'exposant $\ell^{\we_{\si{L}}}$ sur $L$: c'est $L_{\we_{\si{L}}}$.
L'ensemble de cette discussion peut donc se résumer par le diagramme galoisien:\bigskip

\begin{center}
\unitlength=1.5cm
\begin{picture}(6.4,5.2)

\put(0.7,0){$K$}
\put(0.8,0.3){\line(0,1){1.0}}
\put(0.6,1.5){$K_{\we_{\si{K}}}$}
\put(0.8,1.8){\line(0,1){1.0}}
\put(0.6,3.0){$K_{\we_{\si{K}}+\we_{\si{L}}}$}
\put(0.8,3.3){\line(0,1){1.0}}
\put(0.6,4.5){$K_\%$}

\bezier{120}(0.6,0.3)(0.0,2.0)(0.6,4.3)
\put(0.0,2.3){$\Gamma_{\si{\!K}}$}

\put(3.0,1.5){$L=\wH_K$}
\put(3.3,1.8){\line(0,1){1.0}}
\put(3.2,3.0){$L_{\we_{\si{L}}}$}
\put(3.3,3.3){\line(0,1){1.0}}
\put(2.8,4.5){$K^\lc=L_\%$}
\put(1.4,3.05){\line(1,0){1.2}}
\put(1.2,1.55){\line(1,0){1.6}}
\put(1.1,4.55){\line(1,0){1.6}}

\put(1.9,4.65){$\wCl_K$}

\bezier{90}(3.0,1.8)(2.6,3.0)(3.0,4.3)
\put(2.5,3.4){$\Gamma_{\si{\!L}}$}

\put(5.3,3.0){$M=\wH_L$}
\put(5.7,3.3){\line(0,1){1.0}}
\put(5.2,4.5){$L^\lc=M_\%$}
\put(5.35,3.8){$\Gamma_{\si{\!M}}$}

\put(3.7,3.05){\line(1,0){1.5}}
\put(4.4,3.15){$G'$}
\put(3.8,4.55){\line(1,0){1.3}}
\put(4.4,4.65){$\wCl_L$}

\bezier{150}(1.0,0.05)(4, 0.2)(5.6,2.9)
\put(4.0,0.8){$G$}

\end{picture}
\end{center}\medskip\bigskip

Interprétons cela en termes idéliques. Par la théorie du corps de classes, le sous-groupe d'idèles de $\J_K$ qui fixe $L_{\we_{\si{L}}}= LK_{\we_{\si{K}}+\we_{\si{L}}}$ est l'intersection:\smallskip

\centerline{$\J_K^{\ell^{\we_{\si{K}}}}\wU_K\R_K\,\cap\, \J_K^{\ell^{\we_{\si{K}}+\we_{\si{L}}}}\wJ_K \,=\, \J_K^{\ell^{\we_{\si{K}}+\we_{\si{L}}}}(\J_K^{\ell^{\we_{\si{K}}}}\wU_K\R_K\,\cap\,\wJ_K)\,=\,\J_K^{\ell^{\we_{\si{K}}+\we_{\si{L}}}} \wJ_K^{\ell^{\we_{\si{K}}}}\wU_K\R_K$.}\smallskip

Et comme $\ell^{\we_{\si{K}}}$ annule $\,\wCl_K\simeq\wJ_K/\wU_K\R_K$, c'est donc tout simplement: $\J_K^{\ell^{\we_{\si{K}}+\we_{\si{L}}}}\wU_K \R_K$. Il suit:\smallskip

\centerline{$G/G'=\Gal(L_{\we_{\si{L}}}/K)\simeq\J_K/\J_K^{\ell^{\we_{\si{K}}+\we_{\si{L}}}}\wU_K\R_K$.}\smallskip

Par ailleurs, le groupe $G'$ est donné directement comme groupe de classes d'idèles de $L$ par:\smallskip

\centerline{$G'=\Gal(M/L_{\we_{\si{L}}})\simeq\J_L^{\ell^{\we_{\si{L}}}}\wJ_L/\J_L^{\ell^{\we_{\si{L}}}}\wU_L\R_L$.}\smallskip

Cela étant, la trivialité du transfert $G/G' \to G'$ donnée par la théorie des groupes et transportée par le corps de classes nous dit que le morphisme d'extension $j_{\si{L/K}}$ envoie le groupe $\J_K$ dans le sous-groupe $\J_L^{\ell^{\we_{\si{L}}}}\wU_L\R_L$ de $\J_L$, donc son sous-groupe de degré nul $\wJ_K$ dans  $\J_L^{\ell^{\we_{\si{L}}}}\wU_L\R_L\cap\wJ_L$, i.e. dans $\wJ_L^{\ell^{\we_{\si{L}}}}\wU_L\R_L=\wU_L\R_L$, puisque 
$\ell^{\we_{\si{L}}}$ annule $\,\wCl_L\simeq\wJ_L/\wU_K\R_L$ par définition de $\we_{\si{L}}$. Ainsi:

\begin{Th}
Sous la conjecture de Gross-Kuz'min dans $L=\wH_K$, le $\ell$-groupe fini $\,\wCl_K$ des classes logarithmiques (de degré nul) de $K$ capitule dans $\,\wCl_L$.
\end{Th}

\section{Preuve inconditionnelle du Théorème principal}

Venons-en maintenant au cas général formellement identique mais techniquement plus délicat: 
notons $Z_K$ le compositum des $\Zl$-extensions de $K$ contenues dans $K^\lc$ et $Z_L$ son analogue pour $L$; puis $\Z_K$ le sous-groupe d'idèles de $\J_K$ correspondant à $Z_K$, i.e. le noyau du morphisme $\J_K \to \Phi_K=\Gal(Z_K/K)\simeq\Zl^{1+\delta_{\si{K}}}$ donné par le corps de classes et $\Z_L$ son analogue pour $L$, qu'il convient ici de distinguer de $\Z_{L/K}=\Ker (\J_L \to \Phi_{L/K}=\Gal(LZ_K/L)\simeq\ell^{\we_{\si{K}}}\Zl^{1+\delta_{\si{K}}})$.

Observons que le morphisme d'extension $j_{\si{L/K}}$ envoie donc $\Z_K=\sqrt{\wU_K\R_K^\ph}$ dans $\Z_L=\sqrt{\wU_L\R_L^\ph}$.\smallskip

\noindent Désignons enfin par $Z_K^n$ la sous-extension maximale d'exposant $\ell^n$ de $Z_K$ et par  $Z_L^n$ celle de $Z_L$. Partant de $K$, posant $L=\wH_K$ puis $M=\wH_L$ et $G=\Gal(M/K)$, nous obtenons le diagramme:

\begin{center}
\unitlength=1.5cm
\begin{picture}(7.2,6.5)

\put(0.7,0){$K$}
\put(0.8,0.3){\line(0,1){1.0}}
\put(0.6,1.5){$Z_K^{\we_{\si{K}}}$}
\put(0.8,1.8){\line(0,1){2.5}}
\put(0.6,4.5){$Z_K$}

\bezier{120}(0.6,0.3)(0.0,2.0)(0.6,4.3)
\put(0.0,2.3){$\Phi_{\si{\!K}}$}

\put(2.5,1.5){$L=\wH_K$}
\put(2.8,1.8){\line(0,1){1.0}}
\put(2.6,3.0){$LZ_K^{\we_{\si{K}}+\we_{\si{L}}}$}
\put(2.8,3.3){\line(0,1){1.0}}
\put(2.4,4.5){$K^\lc\!=LZ_K$}

\bezier{90}(2.5,1.8)(2.1,3.0)(2.5,4.3)
\put(1.8,3.0){$\Phi_{\si{\!L/K}}$}

\put(1.2,1.55){\line(1,0){1.1}}
\put(1.1,4.55){\line(1,0){1.1}}

\put(1.5,4.7){$\wCl_K^{\,\tor}$}

\put(4.9,4.6){$\Phi_{\si{\!L}}^{\ell^{\we_{\si{L}}}}$}

\put(3.5,3.05){\line(1,0){1.1}}
\put(4.7,3.0){$Z_L^{\we_{\si{L}}}$}
\put(4.8,3.3){\line(0,1){2.5}}
\put(4.7,6.0){$Z_L$}

\put(6.3,3.0){$M=\wH_L$}
\put(6.7,3.3){\line(0,1){2.5}}
\put(6.2,6.0){$L^\lc=MZ_L$}

\put(5.1,3.05){\line(1,0){1.1}}
\put(4.75,2.35){$G'$}
\put(5.1,6.05){\line(1,0){1.0}}
\put(5.5,6.2){$\wCl_L^{\,\tor}$}

\bezier{110}(3.2,2.8)(4.7,2.4)(6.2,2.8)
\bezier{150}(1.0,0.05)(4,0.2)(6.7,2.9)
\put(4.0,0.6){$G$}

\end{picture}
\end{center}\medskip

Reprenons dans ce nouveau cadre les calculs précédents: le groupe d'idèles de $K$ qui fixe le compositum $LZ_K^{\we_{\si{K}}+\we_{\si{L}}}$ est l'intersection de $\J_K^{\ell^{\we_{\si{K}}+\we_{\si{L}}}}\Z_K$ qui fixe $Z_K^{\we_{\si{K}}+\we_{\si{L}}}$ et de $\J_K^{\ell^{\we_{\si{K}}}}\wU_K\R_K$ qui fixe $L$.\par
\noindent D'où:

\centerline{$G/G'=\Gal(LZ_K^{\we_{\si{K}}+\we_{\si{L}}}/K) \simeq \J_K/\J_K^{\ell^{\we_{\si{K}}+\we_{\si{L}}}}(\Z_K\cap\J_K^{\ell^{\we_{\si{K}}}}\wU_K\R_K) = \J_K/\J_K^{\ell^{\we_{\si{K}}+\we_{\si{L}}}}\Z_K^{\ell^{\we_{\si{K}}}}\wU_K\R_K$} 

\noindent i.e.

\centerline{$G/G'\simeq \J_K/\J_K^{\ell^{\we_{\si{K}}+\we_{\si{L}}}}\wU_K\R_K$,}\smallskip

\noindent puisque l'on a $\Z_K^{\ell^{\we_{\si{K}}}} \subset\wU_K\R_K$ par définition de $\we_{\si{K}}$. Par ailleurs, il vient directement:
\smallskip

\centerline{$G'=\Gal(M/LZ_K^{\we_{\si{K}}+\we_{\si{L}}}) \simeq \J_L^{\ell^{\we_{\si{L}}}}\Z_{L/K}/\J_L^{\ell^{\we_{\si{L}}}}\wU_L\R_L$.}\medskip
Cela étant, toujours du fait de la trivialité du transfert $G/G' \to G'$, on conclut que le morphisme d'extension $j_{\si{L/K}}$ envoie le groupe d'idèles $\J_K$ dans le sous-groupe $\J_L^{\ell^{\we_{\si{L}}}}\wU_L\R_L$  de $\J_l$; et, par conséquent, $\Z_K$ dans $\Z_L\cap\J_L^{\ell^{\we_{\si{K}}}}\wU_L\R_L = \Z_L^{\ell^{\we_{\si{L}}}}\wU_L\R_L$; donc finalement dans $\wU_L\R_L$, puisque l'on a $\Z_L^{\ell^{\we_{\si{L}}}} \subset\wU_L\R_L$ par définition de $\we_{\si{L}}$.\smallskip

Il suit, comme précédemment, que $\,\wCl_K^{\,\tor}=\Z_K/\wU_K\R_K$ capitule dans $\,\wCl_L^{\,\tor}=\Z_L/\wU_L\R_L$:

\begin{Th}
Indépendamment de la conjecture de Gross-Kuz'min,, le $\ell$-sous-groupe de torsion $\,\wCl_K^{\,\tor}$ du pro-$\ell$-groupe des classes logarithmiques de $K$ capitule dans $\,\wCl_L^{\,\tor}$ pour $L=\wH_K$.
\end{Th}

\newpage
\Remarque Si $K$ est un corps de nombres totalement réel qui vérifie la conjecture de Leopoldt (par exemple un corps abélien réel), il est bien connu qu'il satisfait aussi la conjecture de Gross-Kuz'min. Et, dans ce cadre, comme établi dans \cite{J60,Ng3}, la conjecture de Greenberg \cite{Gre} revient à
postuler que le $\ell$-groupe des classes logarithmiques $\,\wCl_K$ capitule à un niveau fini $K_n$ de la $\Zl$-extension cyclotomique $K^c=K_\%$. Le Théorème principal de cette note prouve qu'il capitule à un niveau fini de la pro-$\ell$-extension abélienne {\em localement} cyclotomique $K^\lc$.\medskip

\noindent{\em Addendum.} Sous la conjecture de Gross-Kuz'min, il est construit dans \cite{J32} et étudié dans \cite{J36,J39}, pour tout corps de nombres $K$, une $\ell$-tour localement cyclotomique naturelle en prenant $K_1=\wH_K$ puis en itérant le procédé en posant $K_{n+\si{1}}=\wH_{K_n}$ pour $n\ge 1$. La définition de $\wH_K$ donnée dans la présente note permet de construire une telle tour de façon inconditionnelle. Les tours indéfinies au sens de \cite{J32} apparaissent alors comme une catégorie particulière (conjecturalement vide) des tours finies ou infinies ainsi obtenues, suivant que le procédé fait apparaître ou non un étage fini $K_n$ ayant un groupe $\wCl_{K_n}^{\,\tor}$ trivial.

\def\refname{\small{\sc Références}}
{\footnotesize

\bigskip

\noindent{\sc Adresse:}
Univ. Bordeaux \& CNRS,\\
Institut de Mathématiques de Bordeaux,\\
351 Cours de la Libération,\\
F-33405 Talence cedex

\noindent{\sc Courriel:}
 \tt jean-francois.jaulent@math.u-bordeaux.fr
}


\begin{thebibliography}{tt}

\bibitem{Art} {\sc E. Artin}, 
\textit{Idealklassen in Oberkörpern und allgemeine Rezprozitätsgesetz},
Abh. Math. Sem. Hamburg {\bf 7} (1930), 46--51.

\bibitem{A-T} {\sc E. Artin, J. Tate};
\textit{Class Field Theory},
Addison-Wesley (1968).

\bibitem{BJ} {\sc K. Belabas, J.-F. Jaulent},
{\em The logarithmic class group package in PARI/GP},
Pub. Math. Besançon (2016), 5--18.

\bibitem{Fur} {\sc P. Furtwängler},
\textit{Beweis des Hauptidealsatzes für die Klassenkörper algebraischer Zahlkörper},
Abh. Math. Sem. Hamburg {\bf 7} (1930), 14--36.

\bibitem{Gra}{\sc G. Gras}, 
\textit{Class Field Theory: From Theory To Practice},
{Springer-Verlag, 2003}.

\bibitem{Gre}{\sc R. Greenberg}, 
\textit{On a certain $\ell$-adic representation},
Inv. Math. {\bf 21} (1973) 117--124.

\bibitem{Gro}{\sc B. Gross},
{\em $p$-adic $L$-series at $s = 0$},
J. Fac. Sci. Univ. Tokyo IA {\bf 28} (1981), 979--994.

\bibitem{Hil} {\sc D. Hilbert}
{\textit Théorie des corps de nombres algébriques},
Ann. Fac. Sci. Toulouse {\bf 2} (1910), 225--456.

\bibitem{J10} {\sc J.-F. Jaulent}, 
\textit{Sur l'indépendance  $\ell$-adique de nombres algébriques}, 
J. Numb. Th. {\bf 20} (1985), 149--158 .

\bibitem{J28} {\sc J.-F Jaulent},
\textit{Classes logarithmiques des corps de nombres},
J. Théor. Nombres Bordeaux {\bf 10} (1994), 301--325.

\bibitem{J31}{\sc J.-F. Jaulent},
\textit{Théorie $\ell$-adique globale du corps de classes},
J. Théor. Nombres Bordeaux {\bf 10} (1998), 355--397.

\bibitem{J54}{\sc  J.-F. Jaulent},
\textit{Classes logarithmiques et capitulation}, 
Functiones et Approximatio {\bf 54} (2016), 227--239.

\bibitem{J60}{\sc  J.-F. Jaulent},
\textit{Note sur la conjecture de Greenberg},
J. Ramanujan Math. Soc. {\bf 34} (2019), 5--80.

\bibitem{J36} {\sc J.-F. Jaulent, C. Maire,}
{\em A propos de la tour localement cyclotomique d'un corps de nombres,}
Abh. Math. Sem Hamburg {\bf 70} (2000), 239--250.

\bibitem{J39} {\sc J.-F. Jaulent, C. Maire,}
{\em Radical hilbertien et tour localement cyclotomique,}
Japan J. Math. {\bf 28} (2002), 203--213.

\bibitem{J32} {\sc J.-F Jaulent, F. Soriano,}
{\em Sur les tours localement cyclotomiques,}
Archiv der Math. {\bf 73} (1999), 132--140.

\bibitem{Kuz}  {\sc L. V. Kuz'min}, 
{\em The Tate module for algebraic number fields},
 Math. USSR Izv. {\bf 36} (1972), 267--327.
 
 \bibitem{Ng3} {\sc T. Nguyen Quang Do},  
{\em Formule des genres et conjecture de Greenberg},
Annales Math. Québec {\bf 42} (2018), 267--280.

\end{thebibliography}
\end{document}